\documentclass[12pt,leqno,fleqn,epsfig]{article}
\usepackage{amssymb, epsfig, amsmath, amsthm}
\usepackage{mathrsfs}
\usepackage{color}       

\newcommand{\R}{\mathbb{R}}

\newcommand{\N}{\mathbb{N}}

\newcommand{\bc}{\boldsymbol c}
\newcommand{\bd}{\boldsymbol d}
\newcommand{\bfe}{\boldsymbol e}
\newcommand{\bbf}{\boldsymbol f}
\newcommand{\bg}{\boldsymbol g}
\newcommand{\bh}{\boldsymbol h}
\newcommand{\bi}{\boldsymbol i}
\newcommand{\bj}{\boldsymbol j}
\newcommand{\bk}{\boldsymbol k}
\newcommand{\bl}{\boldsymbol l}
\newcommand{\bm}{\boldsymbol m}
\newcommand{\bn}{\boldsymbol n}
\newcommand{\bo}{\boldsymbol o}
\newcommand{\bp}{\boldsymbol p}
\newcommand{\br}{\boldsymbol r}
\newcommand{\bs}{\boldsymbol s}
\newcommand{\bt}{\boldsymbol t}
\newcommand{\bu}{\boldsymbol u}
\newcommand{\bv}{\boldsymbol v}
\newcommand{\bw}{\boldsymbol w} 
\newcommand{\bx}{\boldsymbol x}
\newcommand{\by}{\boldsymbol y}
\newcommand{\bz}{\boldsymbol z}
\newcommand{\buu}{ \underline{\boldsymbol u}}

\newcommand{\bA}{\boldsymbol A}
\newcommand{\bB}{\boldsymbol B}
\newcommand{\bC}{\boldsymbol C}
\newcommand{\bD}{\boldsymbol D}
\newcommand{\bE}{\boldsymbol E}
\newcommand{\bF}{\boldsymbol F}
\newcommand{\bG}{\boldsymbol G}
\newcommand{\bH}{\boldsymbol H}
\newcommand{\bI}{\boldsymbol I}
\newcommand{\bJ}{\boldsymbol J}
\newcommand{\bK}{\boldsymbol K}
\newcommand{\bL}{\boldsymbol L}
\newcommand{\bM}{\boldsymbol M}
\newcommand{\bO}{\boldsymbol O}
\newcommand{\bP}{\boldsymbol P}
\newcommand{\bQ}{\boldsymbol Q}
\newcommand{\bR}{\boldsymbol R}
\newcommand{\bS}{\boldsymbol S}
\newcommand{\bT}{\boldsymbol T}
\newcommand{\bU}{\boldsymbol U}
\newcommand{\bV}{\boldsymbol V}
\newcommand{\bW}{\boldsymbol W}
\newcommand{\bX}{\boldsymbol X}
\newcommand{\bY}{\boldsymbol Y}
\newcommand{\bZ}{\boldsymbol Z}

\newcommand{\mes}{\operatorname{\rm meas}}    
\newcommand{\esssup}{\operatorname*{ess\,sup}}

\newcommand{\osc}{\operatorname*{osc}}

\newcommand{\supp}{\operatorname*{supp}}

\newcommand{\const}{\operatorname*{const}}

\newcommand{\bb}{\begin{equation}}
\newcommand{\ee}{\end{equation}}
\newcommand{\bq}{\begin{eqnarray}}
\newcommand{\eq}{\end{eqnarray}}
\newcommand{\bqn}{\begin{eqnarray*}}
\newcommand{\eqn}{\end{eqnarray*}}

\newcommand{\var}{\varepsilon}

\newcommand{\intl}{\int\limits}

\newcommand{\Beweisende}{\rule{0.2cm}{0.2cm}}

\newcounter{secnum}

\newtheorem{thm}{Theorem}[section]
\newtheorem{cor}[thm]{Corollary}
\newtheorem{lem}[thm]{Lemma}

\theoremstyle{definition}

\newtheorem{rem}[thm]{Remark}

\title{On the local  Type I conditions for \\ the 3D Euler equations} 
 
\author{Dongho Chae$^*$  and J\"{o}rg Wolf $^\dagger$\\
\ \\
Department of Mathematics\\
Chung-Ang University\\
 Seoul 156-756, Republic of Korea\\
 ($*$)e-mail: dchae@cau.ac.kr\\
($\dagger$)e-mail: jwolf2603@cau.ac.kr}
\date{}

\begin{document}
\maketitle
\begin{abstract}
We prove local  non blow-up theorems  for the  3D incompressible Euler equations under local Type I conditions. More specifically, for a classical solution $v\in L^\infty (-1,0; L^2 ( B(x_0,r)))\cap L^\infty_{\rm loc} (-1,0; W^{1, \infty} (B(x_0, r)))$ of the 3D Euler equations, where $B(x_0,r)$ is the ball with radius $r$ and the center at $x_0$, if the  limiting values of certain  scale invariant quantities for a solution $v(\cdot, t)$ as $t\to 0$ are small enough, then  $ \nabla v(\cdot,t) $ does not blow-up at $t=0$ in $B(x_0, r)$.
\\
\
  \noindent{\bf
keywords:} Euler equation, finite time blow-up, energy concentration, discretely self-similar solution

\end{abstract}

\section{Introduction}
\label{sec:-1}
\setcounter{secnum}{\value{section} \setcounter{equation}{0}
\renewcommand{\theequation}{\mbox{\arabic{secnum}.\arabic{equation}}}}

We consider the 3D homogeneous  incompressible Euler equation in a cylinder $ Q= \R^{3}\times (-1,0)$
\begin{align}
\partial _t v + v\cdot \nabla v &= - \nabla p \quad  \text{ in}\quad  Q, 
\label{eul1}\\
\nabla \cdot v &=0 \quad  \text{ in}\quad  Q, 
\label{eul2}
\end{align}
where $ v=(v_1(x,t), v_2(x,t), v_3(x,t))$ stands for the velocity of the fluid and $ p=p(x,t)$ stands for the pressure. 
The local in time well-posedness in  the Sobolev space $W^{k,p} (\Bbb R^3)$, $k>3/p+1$, $1<p<+\infty$, for the Cauchy problem of the system \eqref{eul1}-\eqref{eul2}   is well-known due to the proof by Kato-Ponce\cite{kat}. The question of the spontaneous apparition of singularity 
from the local in time smooth solution, however, is still an outstanding open problem in the mathematical fluid mechanics(see e.g.\cite{maj, con1} for surveys of studies devoted to the problem).  We say a local in time smooth solution $v\in C([-1, 0); W^{k,p}(\Bbb R^3))$, $k>3/p+1$, $1<p<+\infty$,
does not blow up (or becomes regular) at $t=0$ if 
\begin{equation}\label{noblow}
\limsup_{t\to 0^-} \|v(t)\|_{W^{k,p}(\Bbb R^3)}<+\infty.
\end{equation}
It is easy to show from the local in time well-posedness estimates that \eqref{noblow} is guaranteed if
\begin{equation}\label{noblow1}
\int_{-1} ^0 \|\nabla v(t)\|_{L^\infty(\Bbb R^3)} dt <+\infty.
\end{equation} 
The celebrated Beale-Kato-Majda criterion\cite{bea} shows that one can replace \eqref{noblow1} by a weaker condition
\begin{equation}\label{noblow2}
\int_{-1} ^0 \|\omega (t)\|_{L^\infty(\Bbb R^3)} dt <+\infty, \qquad \omega =\nabla \times v.
\end{equation} 
(see also \cite{con2, den} for geometric type criterion, and \cite{koz} for a refinement of \eqref{noblow2}).
The conditions \eqref{noblow1} or \eqref{noblow2} can be regarded as regularity criteria of the Serrin type in the Navier-Stokes equations.
There exist also another form of {\em local}  regularity criteria, called $\varepsilon-$regularity criteria, which claims that if certain scaling invariant quantities are
small enough in a local space-time neighborhood, then weak solution becomes regular in the neighborhood. A typical example of such smallness condition, introduced  by Caffarelli, Kohn and Nirenberg in \cite{caf}, which guarantees the regularity  near $(x,t)=(0,0)$ for a {\em suitable weak solution}  of the Navier-Stokes equations  is
\begin{equation}\label{ckn}
\limsup_{r \to 0^+} \frac{1}{r} \int_{-r^2 }^0 \int_{\{|x|<r\}} |\nabla v(x,t)|^2dxdt <\varepsilon,
\end{equation} 
where $\varepsilon >0$ is an absolute constant.
The replacement of $\varepsilon$ by  finite constant $C$ in \eqref{ckn} is called local Type I condition for the Navier-Stokes equations
 (cf. \cite{tsai3, ser}).
In  view of the scaling property of the Euler equations a natural local Type I condition with smallnes, which guarantee no blow-up at $t=0$
for a classical solution $v\in C([-1, 0); W^{1, \infty } (\Bbb R^3))$ would be
\begin{equation}\label{ckn1}
\limsup_{t\to 0^-} (-t)\|\nabla v(t)\|_{L^\infty (B(r))}<\varepsilon,
\end{equation} 
where we used $B(r)=B(0,r)$ with   $B(x_0, r)=\{ x\in \Bbb R^3\, |\, |x-x_0|<r\} $.  
Indeed in \cite{cha2}(see also \cite{hou} for an independent result) it has been shown that if 
\begin{equation}\label{ckn2}
\limsup_{t\to 0^-} (-t)\|\nabla v(t)\|_{L^\infty (\Bbb R^3)}< 1,
\end{equation} 
then  there exists no blow-up at $t=0$ for a classical solution to the Euler equations on $\Bbb R^3\times (-1,0)$.
Our first aim in this paper is to {\em localize in space}  \eqref{ckn2}, and prove the following theorem.

\begin{thm}
\label{thm1.1}
Let $ v \in L^\infty(-1, 0; L^2( B(r)))\cap L^\infty_{ loc}([-1,0); W^{1,\infty}(B(r)))$ be a solution to the Euler equations 
\eqref{eul1}, \eqref{eul2} with $ v(-1)\in W^{2,\, p_0}(B(r))$ for some $ 3 < p_0<+\infty$.  
We assume there exists $r_0\in (0, r)$ such that 
\begin{equation}
\limsup_{ t\to  0^-}(-t)\| \nabla v(t)\|_{ L^\infty(B(r_0))}  <1. 
\label{decay}
\end{equation}  
Then $\lim\sup_{t\to 0^-} \| v( t)\|_{ W^{2,\, p_0}(B(\rho))}<+\infty$ for all $ \rho\in (0, r_0)$. 
\end{thm}

 The proof of the above theorem is given in the Section 2 to Section 4.
From the structure of the Euler equations the estimation of  the $ L^p_{\rm loc} (\Bbb R^3)$ norm of second  derivatives usually are  
obtained by means of Gronwall's Lemma. In oder to handle the integrals involving derivatives with cut off function
it was crucially helpful  to introduce the following  transformation of the solutions $v(x,t)  \mapsto w(y,t)= v(  (1+ (-t)^\theta ) y,t)$ for appropriately chosen $ 0<\theta <1$. 
 
\vspace{0.2cm}
\hspace{0.5cm}
In our  second main result below  we use  Theorem \ref{thm1.1} to deduce that local small oscillation near $t=0$ implies also no blow-up of  a classical solution on $ B(r)\times (-1,0)$.
\begin{thm}
\label{thm1.2}
Let $ v \in L^\infty(-1, 0; L^2(B(r)))\cap L^\infty_{ loc}([-1,0); W^{1,\infty}(B(r)))$ be a solution to the Euler equations 
\eqref{eul1}, \eqref{eul2}  with $ v(-1)\in W^{2,\, p_0}(B(r))$ for some $ 3 < p_0<+\infty$. We assume there exists $ 0< r_0 <  r/2$ such that 
\begin{equation}
\limsup_{ t\to0^-} (-t)\sup_{ x_0\in B(r/2)} \osc_{ B(x_0, r_0 (-t)^{ 2/5} )} (\nabla v(t))    <1. 
\label{odecay}
\end{equation}  
Then $\limsup_{t\to 0^-} \| v(t)\|_{ W^{2,\, p_0}(B(r/4))} <+\infty$.
\end{thm} 

The key ingredient in the proof of Theorem\,\ref{thm1.2}  is the fact that under Type I condition(replacing $\varepsilon $ by any finite constant $C$ in \eqref{ckn1}) there exists no atomic energy concentration in $B(r)$,  which is proved in \cite{cw2}, that makes  the local energy $\int_{B(x_0,r)}|v(x,t)|^2 dx$ uniformly  small with respect to $(x_0, t)\in \Bbb R^3 \times (-1,0)$ for small $r>0$.

\begin{rem}
Below we present two sufficient conditions on  $v$,  which imply \eqref{odecay}. The first one is obvious.  If there exists 
a function $ \eta : (0, \infty) \rightarrow \R$ with  $ \eta (r) \rightarrow 0$ as $ r \rightarrow 0$ such that
\begin{equation}
(-t)| \nabla v(x, t)- \nabla v(y, t) | \le \eta\Big( \frac{| x-y|}{(-t)^{ 2/5}}\Big)\quad \forall\,x,y\in \R^{3}, t\in (-1,0),  
\label{decay1}
\end{equation}
then, $v$ satisfies  \eqref{odecay}.  
The second condition is given in terms of the Fourier transform.  Given $\delta \in (0, 1)$, if  there exists $ 0<R_0<+\infty$ such that 
\begin{equation}
\intl_{ \R^{3}  \setminus B(R_0 (-t)^{ -2/5})} | \xi | | \mathcal{F} v(\xi , t)| d\xi  \le  \frac{1- \delta }{-2t}  \quad    \forall\,t\in (-1,0), 
\label{decay2}
\end{equation} 
then  the condition  \eqref{odecay} for $v$ follows.  Indeed, let $ g(\xi , t) = \mathcal{F} (\nabla v(\xi ))$, then we see that 
\begin{align*}
&\nabla v(x,t)- \nabla v(y,t) 
\\
&= \mathcal{F}^{ -1}g(\cdot ,t)(x)- \mathcal{F}^{ -1}g(\cdot ,t)(y)
\\
&= \intl_{ \R^{3}  \setminus  B(R_0 (-t)^{ -2/5})} (e^{ 2\pi  i x\cdot  \xi } - e^{ 2\pi i   y\cdot  \xi }) g(\xi,t ) d\xi  + \intl_{ B(R_0 (-t)^{ -2/5})} (e^{ 2 \pi i x\cdot  \xi } - e^{ 2\pi  i y\cdot  \xi }) g(\xi,t ) g(\xi,t ) d\xi  
\end{align*} 
which leads to the inequality 
\begin{equation}
| \nabla v(x,t)- \nabla v(y,t)| \le  \frac{1- \delta}{-t}  + \intl_{ B(R_0 (-t)^{ -2/5})} | e^{2 \pi  i(x-y)\cdot  \xi } - 1| | g(\xi,t )| d\xi.   
\label{decay3}
\end{equation}
For $ | x-y| \le r_0 (-t)^{ 2/5}$ and  $ 0< r_0 \le R_0^{ -1}$ the second term can be estimated as follows\begin{align*}
\intl_{ B(R_0 (-t)^{ -2/5})} | e^{2 \pi i (x-y)\cdot  \xi } - 1| | g(\xi,t )| d\xi 
&\le 2\pi R_0\frac{| x-y|}{(-t)^{ 2/5}}  \intl_{ \R^{3}}  | \xi || \mathcal{F}v( \xi , t )| d\xi  
\\
&\le 2\pi r_0 R_0 \intl_{ \R^{3}}  | \xi || \mathcal{F}v( \xi , t )| d\xi.
\end{align*}
Thus  if we choose  $ r_0$ such that 
\[
2\pi r_0 R_0 \sup_{ t\in (-1,0)}(-t)\intl_{ \R^{3}}  | \xi || \mathcal{F}v( \xi , t )| d\xi \le \frac{\delta }{2},
\]  
the condition \eqref{decay} holds with $ \frac{\delta }{2}$ in place of $ \delta $.
Here,  we have used the fact that there exists a constant $ C_1 >0$ such that 
\[
(-t) \intl_{ B(R_0 (-t)^{ -2/5})}  | \xi || \mathcal{F}v( \xi , t )| d\xi \le C_1.
\]
This can be checked  by  H\"{o}lder's inequality and Plancherel's theorem as follows.
\begin{align*}
(-t) \intl_{ B(R_0 (-t)^{ -2/5})}  | \xi || \mathcal{F}v( \xi , t )| d\xi
\le cR_0^{ 5/2} \| \mathcal{F}v(t)\|_{ L^2} = cR_0 ^{ 5/2} E^{ 1/2}.  
\end{align*}

\end{rem}

\section{Uniform smallness of the local energy}
\label{sec:-2}
\setcounter{secnum}{\value{section} \setcounter{equation}{0}
\renewcommand{\theequation}{\mbox{\arabic{secnum}.\arabic{equation}}}}

Our aim in this section is to prove the following result, which is interesting itself.

\begin{thm}
\label{thm2.1}
Let $ v \in L^\infty(-1, 0; L^2(B(1)))\cap L^\infty_{ loc}([-1,0); W^{1,\, \infty}(B(1)))$  be a solution to the Euler equations 
\eqref{eul1}, \eqref{eul2}, which satisfies the following  condition 
\begin{equation}
\sup_{ t\in (-1,0)}(-t) \| \nabla v(t)\|_{L^\infty(B(1))} \le C_0. 
\label{type1}
\end{equation}  
Then for every $ \var >0$ there exists $ 0< \widetilde{R} =\widetilde{R} (\var ) \le  \frac{1}{2}$  such that for all $ y_0\in B(1/2)$ 
it holds 
\begin{equation}
\sup_{ t\in ( -\widetilde{R} ^{ 5/2} , 0)} \intl_{B(y_0, \widetilde{R} )} | v(t)|^2 dx\le \var.  
\label{morrey}
\end{equation}
\end{thm}

The proof of Theorem\,\ref{thm2.1} will be achieved after proving several lemmas.   Given $z_0=(x_0,t_0)\in \Bbb R^3\times (-\infty, 0]$, $ r>0$, we denote $Q(z_0, r)= B(x_0,r)\times (t_0-r^{5/2}, t_0)$, and $Q(r)= Q(0,r)$.
For   $\Omega \subset \Bbb R^3$  by $ W^{\theta ,\, p}(\Omega )$ we denote  the usual Sobolev-Slobodecki\u{i} space, which 
consists of all functions $ f\in L^p(\Omega )$ such that the following semi norm is finite
\[
| f|^p_{ W^{\theta, p}(\Omega )} = \intl_{\Omega} \intl_{\Omega}  \frac{| f(x)- f(y)|^p }{| x-y|^{ 3+p\theta }} dxdy,\quad  
f\in W^{\theta, p}(\Omega ),\quad 0<\theta <1, \,\,p \ge 1. 
\]
\begin{lem}
\label{lem2.2}
Let the   assumption of Theorem\,\ref{thm2.1} be satisfied. Let  $ x_0 \in B(1/2)$. Then,  for every $ \var >0$ 
there exists $ 0< R_0= R_0(x_0, \var )< 1$ such that for all $ 0< R \le  R_0$ 
it holds 
\begin{equation}
R^{ -1} \intl_{Q((x_0,0), R)} | v|^3 dxdt\le \var.  
\label{morreyL3}
\end{equation}

\end{lem}

{\bf Proof}:  We prove the assertion of the theorem by an indirect argument.  To this end let us assume the assertion is false. 
Then there exist $ x_0\in B(1/2)$ and a sequence  $ \{ r_k\}$ of numbers in $ (0,1/2)$,  which converges to zero as $ k \rightarrow +\infty$, satisfying 
\begin{equation}
r_k^{ -1} \intl_{Q((x_0,0), r_k)} | v|^3 dxdt> \var \qquad  \forall\,k\in \N. 
\label{false}
\end{equation}
Without the  loss of generality we may assume $ x_0=0$. 

\hspace{0.5cm}
Since the solution is defined locally, we cannnot  expect global bounds on the pressure. 
By this reason the compactness lemma of Lions-Aubin  type does not work in this situation. This forces us to work with the notion of  local pressure. 
As in \cite{wol, cw2} we introduce the projection $ E_{ B(1)}^{ \ast} : W^{-1,\, q}(B(1)) \rightarrow W^{-1,\, q}(B(1))$ onto the space of functionals given by a gradient $ \nabla \pi $. In fact, here 
 $ \pi \in L^q_0(B(1))$  denotes the pressure of the solution to the Stokes equations in $ B(1)$ with zero boundary data and force $ f$ . 
 We define
 \begin{align*}
 &\nabla p_h = -E^{ \ast}_{ B(1)}(v),\quad  \nabla p_0 = - E^{ \ast}_{ B(1)} ((v\cdot \nabla) v), 
 \\
 &\qquad \qquad \widetilde{v} = v + \nabla p_h.  
 \end{align*}
 From this definition we find easily  that $ (\widetilde{v}, \nabla p_h, p_0 )$ solves the following system in $ B(1)\times (-1,0)$ in the sense of distributions
 \begin{equation}
 \partial _t \widetilde{v} + (v\cdot \nabla )\widetilde{v} - v \cdot \nabla^2 p_h = - \nabla p_0.
 \label{eulerloc}
 \end{equation} 
 
Let $ 0< \rho < +\infty$ be arbitrarily chosen. Recalling that $ \nabla p_h$ is harmonic, by using the mean value property of harmonic functions, we get for all $ t\in (-1,0)$ and for all $ k\in \N$ such that $ \rho r_k \le \frac{1}{2}$ the following estimates
 \begin{align}
\intl_{B(\rho r_k)} | \nabla p_h(t) |^3 dx  &\le c \rho ^3 r_k ^3 \| \nabla p_h(t)\|^3 _{ L^\infty(B(1/2))} \le  
c \rho ^3 r_k ^3 \| \nabla p_h(t) \|^3_{ L^\infty(B(1))} 
 \cr
 &\le c \rho ^3 r_k ^3 \| v \|^3_{ L^\infty(-1,0; L^2(B(1))}. 
 \label{esth}
 \end{align}
Next, we define sequences of scaled velocities and pressures,
\begin{align*}
v_k (x,t)  &= r_{ k}^{ 3/2} v(r_k x, r^{ 5/2} t), 
\\
p_{ h,k} (x,t)  &= r_{ k}^{ 1/2} p_h(r_k x, r^{ 5/2} t), 
\\
p_{ 0,k} (x,t) &= r_{ k}^{ 3} p_0(r_k x, r^{ 5/2} t), \quad  (x,t)\in B(r_k^{ -1})\times (-1,0). 
\end{align*}
Using the transformation formula of the Lebesgue integral, the condition  \eqref{false} becomes   
\begin{equation}
 \intl_{Q( 1)} | v_k |^3 dxdt \ge  \var. 
\label{false1}
\end{equation}

Setting $ \widetilde{v} _k = v_k + \nabla p_{ h,k}$ in  \eqref{eulerloc}, we see that the following equations  are satisfied in the sense of distributions.
 \begin{equation}
 \partial _t \widetilde{v}_k + (v_k\cdot \nabla )\widetilde{v}_k - v_k  \cdot \nabla^2 p_{ h,k} = - 
 \nabla p_{ 0,k} \qquad \text{ in  $\quad B(r_k^{ -1})\times (-1,0)$.}
  \label{eulerloc1}
 \end{equation} 
 From \eqref{esth} we immediately get for all $ 0< \rho < +\infty$ 
and for all $ k\in \N$ such that $ \rho r_k \le  \frac{1}{2}$
\begin{equation}
\| \nabla p_{ h,k}\|_{ L^3(B(\rho )\times (-1,0))}^3 = r_k^{-1} \| \nabla p_{ h}\|_{ L^3(B(\rho r_k)\times (- r_k^{ 5/2},0))}^3
\le c \rho^3 r_{ k}^2. 
\label{esth1}
\end{equation} 
This yields 
\begin{equation}
 \nabla p_{ h,k} \rightarrow 0  \quad  \text{{\it in}}\quad  L^3(B(\rho )\times (-1,0))\quad  \text{{\it as}}\quad  k \rightarrow +\infty
 \quad  \forall\,0< \rho < +\infty.  
\label{limhp}
\end{equation}

Since both the $ L^\infty(-1,0; L^2(B(1)))$ norm and the Type I condition \eqref{type1} are invariant under the above scaling,   we obtain 
for all $ 0< \rho < +\infty$ and for all $ k\in \N$ such that $ \rho r_k \le \frac{1}{2}$
\begin{align}
&\| v_k\|_{ L^\infty(-1,0; L^2(B(\rho )))} \le \| v\|_{ L^\infty(-1,0; L^2(B(1)))}, 
\label{2.1}
\\
&\sup_{ t\in (-1,0)}(-t) \| \nabla v_k(t)\|_{L^\infty(B(\rho ))} \le C_0. 
\label{2.2}
\end{align}

By interpolation between the two bounds \eqref{2.1} and \eqref{2.2}  we see that for all $ 0<\rho < +\infty$ the sequence 
$ \{ v_k\}$ is bounded in 
$ L^3(-1,0; $ $W^{\theta ,\, 3}(B(\rho )))$ for each $ 0 \le  \theta < \frac{1}{3}$.  Using the regularity properties of harmonic functions,  we also see that $ \{ \widetilde{v} _k\}$ is bounded in $ L^3(-1,0; $ $W^{\theta ,\, 3}(B(\rho )))$.  
In particular, $ \{ v_k\}$ is bounded in $ L^3(B(\rho )\times (-1,0))$ which shows that $ \{ p_{ 0, k}\}$ 
is bounded in $ L^{ 3/2}(B(\rho )\times (-1,0))$.  Accordingly, $ \{ \partial _t \widetilde{v}_k \}$ is bounded in 
$ L^{ 3/2}(-1,0; W^{-1,\, 3/2}(B(\rho )))$ for all $ 0< \rho < +\infty$. 
Using Banach-Alaoglu's theorem and applying a compactness lemma due to Simon\cite{sim}, and Cantor's diagonalization principle,   eventually passing to a subsequence, 
we get a  limit $ v^{ \ast} \in L^\infty(-1, 0; L^2(\R^{3}))\cap L^\infty_{ loc}([-1,0); W^{1,\, \infty}(\R^{3})) $ such that 
for all $ 0< \rho < +\infty$
\begin{align}
 \widetilde{v} _k &\rightarrow   v^{ \ast}\quad  \text{{\it weakly-$ \ast$ in}}\quad  L^\infty(-1,0; L^2(B(\rho )))\quad  \text{{\it as}}\quad  k \rightarrow +\infty, 
 \label{lim1}
\\
 \widetilde{v} _k &\rightarrow   v^{ \ast} \quad  \text{{\it in}}\quad  L^3(B(\rho )\times (-1,0))\quad  \text{{\it as}}\quad  k \rightarrow +\infty.
\label{lim2}
\end{align}
Using \eqref{lim2} and \eqref{limhp}, we may  let  $ k \rightarrow +\infty$ in \eqref{false1}, which yields 
\begin{equation}
 \intl_{Q(1)} | v^{ \ast}|^3 dxdt \ge  \var.
\label{false2}
\end{equation}
Furthermore, after passing $ k \rightarrow +\infty$ in \eqref{eulerloc1}, we deduce that $ v^{ \ast} 
\in L^2(-1,0; L^2(\R^{3}))\cap L^\infty_{ loc}([-1,0); W^{1,\, \infty}(\R^{3}))$ is a solution to the Euler equations in $ \R^{3}\times (-1,0)$.  Thanks to \eqref{2.2} we see that the function $ v^{ \ast}$ enjoys the Type I blow up condition 
\begin{equation}
\sup_{ t\in (-1,0)}(-t) \| \nabla v^{ \ast}(t)\|_{L^\infty( \R^{3})} \le C_0. 
\label{2.3}
\end{equation}

Observing \eqref{eulerloc1}, and  noting that $ v_k(t)$ is bounded and Lipschitz before the blow up time we get  the following local energy equality which holds 
for every $ -1 \le t < s < 0$ and for all $ \phi \in C^{\infty}_c(B(\rho ))$ with $ r_k \rho < 1$
\begin{align}
\intl_{ B(\rho )} | \widetilde{v}_k (t)|^2 \phi dx &= \intl_{ B(\rho )} | \widetilde{v}_k (s)|^2 \phi dx
+  \intl_{t}^{s} \intl_{ B(\rho )} (| \widetilde{v} _k|^2 v_k + 2 p_{ 0,k} \widetilde{v}_k)  \cdot \nabla \phi  dxd\tau
\cr
&\qquad \qquad + \intl_{t}^{s} \intl_{B(\rho )} \widetilde{v} _k  \nabla \phi  dx d\tau . 
\label{locen}
\end{align}

In the  discussion below $ \mathcal{M}(\R^{3})$ denotes the space of Radon measures, while $ \mathcal{M}^+(\R^{3})$ stands for the space of positive Radon measures both on $ \R^{3}$. 
As we have proved in \cite{cw2} there exists a unique measure $ \widetilde{\sigma}  \in \mathcal{M}^+(\R^{3})$ such that 
\begin{equation}
 | \widetilde{v} (t)|^2 dx \rightarrow \widetilde{\sigma}    \quad  \text{{\it weakly-$ \ast$ in}}\quad  \mathcal{M}(B(1))\quad  \text{{\it as}}\quad  
 t \rightarrow 0^-. 
\label{2.4}
\end{equation}

This implies 
\[
 | \widetilde{v} _k(s)|^2 dx \rightarrow \widetilde{\sigma} _k   \quad  \text{{\it weakly-$ \ast$ in}}\quad  \mathcal{M}(B(\rho ))\quad  \text{{\it as}}\quad  
 s \rightarrow 0^-,
\]
 where $ \widetilde{\sigma} _k$ is defined as 
\[
\intl_{B(\rho )} \phi d\widetilde{\sigma} _k = \intl_{ B(r_k \rho )} \phi  \Big(\frac{x}{r_k}\Big) d\widetilde{\sigma}.    
\]
Thus, in \eqref{locen} letting $ s \rightarrow 0$, we arrive at 
\begin{align}
\intl_{B(\rho )} | \widetilde{v}_k (t)|^2 \phi dx &= \intl_{B(\rho )} \phi d\widetilde{\sigma} _k
+  \intl_{t}^{0} \intl_{ B(\rho )} (| \widetilde{v} _k|^2 v_k + 2 p_{ 0,k} \widetilde{v}_k)  \cdot \nabla \phi  dxd\tau
\cr
&\qquad \qquad + \intl_{t}^{0} \intl_{ B(\rho )} \widetilde{v} _k  \nabla \phi  dx d\tau . 
\label{locen1}
\end{align}
On the other hand, it can be checked easily that $ \| \widetilde{\sigma} _k\| \le \| \widetilde{\sigma }\|$. Hence, eventually passing to a subsequence,  we get 
 $ \sigma ^{ \ast} \in 
 \mathcal{M}^+(\R^{3})$ such that for all $ 0< \rho <+\infty$
 \[
\widetilde{ \sigma} _k\rightarrow \sigma^{ \ast}  \quad  \text{{\it weakly-$ \ast$ in}}\quad 
  \mathcal{M}(B(\rho ))\quad  \text{{\it as}}\quad  
 s \rightarrow 0^-. 
\]
Our aim is to show that $ \sigma ^{ \ast} = \widetilde{\sigma} (\{ 0\}) \delta _0$.  Arguing as in \cite{cw2}, we infer 
for $ \phi \in C^{\infty}_c(\R^{3})$ with $ \supp(\phi ) \subset B(\rho )$
\begin{align*}
\intl_{ \R^{3}} \phi d\sigma ^{ \ast} &=  \lim_{k \to \infty}\intl_{B(\rho )} \phi d\widetilde{\sigma} _k =
 \lim_{k \to \infty}\intl_{ B(r_k \rho )} \phi \Big(\frac{x}{r_k}\Big) d\widetilde{\sigma}  
 \\
 & =\lim_{k \to \infty}\intl_{B( r_k\rho )  \setminus \{ 0\}} \Big\{\phi \Big(\frac{x}{r_k}\Big) - \phi (0)\Big\}d\widetilde{\sigma}  +  
 \phi (0)\lim_{k \to \infty} \intl_{B( r_k\rho )} d\widetilde{\sigma} 
 \\
  & = \phi (0) \widetilde{\sigma}  (\{ 0\}).   
\end{align*}
This shows the claim. Thus, we are in a position to apply \cite[Theorem\,3.1]{ cw2},  which excludes the concentration of energy  at one point. 
Hence, $ v^{ \ast} $ must vanish. However this contradicts with \eqref{false2}. Accordingly the assertion must be true.    \hfill \Beweisende 

\vspace{0.5cm}  
Next, we show the smallness of the local energy. 

\begin{lem}
\label{lem2.3}
Let all assumptions of Theorem\,\ref{thm2.1} be fulfilled. 
Let $ x_0\in B(1/2)$. Then for every $ \var >0$ there exists $ 0< R_1 =  R_1(\var, x_0 )< \frac{1}{2}$ such that  
\begin{equation}
\| v(t)\|_{ L^2(B(x_0, R_1))} \le \var \quad  \forall\, - R_1^{ 5/2} \le t <0. 
\label{energy}
\end{equation} 
\end{lem}

{\bf Proof}: Let $ 0< \delta <1$ be a number,  which will be specified below. By Lemma\,\ref{lem2.2}  there exists 
$ 0< R_0=R_0((\delta \var)^{ 3/2}, x_0 ) <1$ such that \eqref{morreyL3} holds with $ (\delta \var)^{ 3/2}$ in place of $ \var $. Accordingly, by the help of H\"{o}ler's inequality   we find 
\begin{equation}
R_0^{ -5/2} \intl_{Q((x_0,0), R_0)} | v|^2 dxdt  \le | B(1)|^{ 1/3} \delta \var.  
\label{morrey1}
\end{equation}
Thus, thanks to the mean value property of the integral we may choose $ s_0 \in (-R_0 ^{5/2}, 0)$ such that 
\begin{equation}
 \intl_{B(x_0, R_0)} | v(s_0)|^2 dx \le | B(1)|^{ 1/3} \delta \var.  
\label{morrey2}
\end{equation}

As in the proof of Lemma\,\ref{lem2.2} we define  the local pressure 
\[
\nabla p_{ h} = - E^{ \ast}_{ B(x_0, R_0)}(v),\quad \nabla p_{ 0} = - E^{ \ast}_{ B(x_0, R_0)}((v\cdot \nabla) v),
\] 
and set $ \widetilde{v} = v+ \nabla p_h$. As above we see that the function $ \widetilde{v} $ solves the modified 
Euler equations 
\begin{equation}
\partial _t \widetilde{v} + v\cdot \nabla \widetilde{v}  -v\cdot \nabla^2 p_h = - \nabla p_0\quad  \text{ in} \quad  
B(x_0, R_0)\times (-1,0)
\label{eul3}
\end{equation}
in the sense of distributions. Furthermore the following local energy identity holds true for all $ -1 \le t <s <0$ and for 
all $ \phi \in C^{\infty}_c(B(x_0, R_0))$
\begin{align}
&\intl_{ B(x_0, R_0)} | \widetilde{v} (t)|^2 \phi dx = \intl_{ B(x_0, R_0)} | \widetilde{v} (s)|^2 \phi dx
+  \intl_{t}^{s} \intl_{ B(x_0, R_0)} (| \widetilde{v} |^2 v + 2 p_0 v) \cdot \nabla \phi  dxd\tau
\cr
&\qquad \qquad +  \intl_{t}^{s} \intl_{ B(x_0, R_0)} v \cdot \nabla ^2 p_h\cdot \widetilde{v}   \phi  dxd\tau. 
\label{locen2}
\end{align}
Now, let $ \phi \in C_c ^\infty (\Bbb R^3)$ be a radial cut-off function  such that $ \phi \equiv 1$ 
on $ B(x_0, R_0/2)$, $ 0 \le \phi \le 1$ in $ B(x_0, R_0)$ and $ | \nabla \phi | \le c R_0^{ -1}$.  We set $ s=s_0$ 
and $ t\in (-R_0^{ 5/2},0)$ in \eqref{locen2}. Using \eqref{morrey2}, we see that 
\begin{align}
&\intl_{ B(x_0, R_0)} | \widetilde{v} (t)|^2 \phi dx \le  | B(1)|^{ 1/3} \delta \var
+ c R_0^{ -1} \intl_{ Q((x_0,0), R_0)} (| \widetilde{v} |^2 | v| + 2 | p_0| | v| )  dxd\tau
\cr
&\qquad \qquad +  \intl_{t}^{s_0} \intl_{ B((x_0, 0), R_0)} | v|| \nabla ^2 p_h| | \widetilde{v}|   \phi  dxd\tau
= | B(1)|^{ 1/3} \delta \var + I+II. 
\label{locen3}
\end{align}   
To estimate the first integral we make use of the pressure estimates
\begin{align*}
\| \nabla p_h(t)\|_{ L^3(B(x_0, R_0))} &\le c \|v(t)\|_{ L^3(B(x_0, R_0))},\quad  
\\
\|  p_0(t)\|_{ L^{ 3/2}(B(x_0, R_0))} &\le c \| v(t)\|^2_{ L^3(B(x_0, R_0))}, 
\end{align*}
 which together with \eqref{morrey1} shows that 
 \[
I \le c R^{ -1}_0  \intl_{ Q((x_0,0), R_0)} | v|^3   dxd\tau \le c(\delta \var )^{ 3/2}. 
\]
To estimate the second integral  we first apply H\"older's inequality to get 
\begin{align*}
II \le \bigg(\intl_{ Q((x_0,0), R_0)} | v|^3   dxd\tau\bigg)^{ 2/3}  
\bigg(\intl_{ Q((x_0,0), R_0)} | \nabla^2 p_h|^3 \phi ^3   dxd\tau\bigg)^{ 1/3}.
\end{align*}
Applying Sobolev's embedding theorem together with \cite[Lemma\,A.1]{cw2}, 
we estimate for $  \tau \in (-R_0^{5/2},0)$
\begin{align*}
&\| \nabla^2 p_h(\tau ) \phi \|_{ L^3(B(x_0, R_0))}   \\
& \le c R_0^{ -1/2} \| \nabla^2 p_h(\tau ) \phi \|_{ L^2(B(x_0, R_0))} 
+ c \| \nabla^2 p_h(\tau ) \nabla \phi \|_{ L^2(B(x_0, R_0))}
\\
&\qquad \qquad \qquad + c \| \nabla^3 p_h(\tau ) \phi \|_{ L^2(B(x_0, R_0))}
\\
& \le c R_0^{ -3/2} \| \nabla p_h(\tau ) \|_{ L^2(B(x_0, R_0))} 
\\
& \le c R_0^{ -3/2} \| v(\tau ) \|_{ L^2(B(x_0, R_0))}  \le  c R_0^{ -1} \| v(\tau ) \|_{ L^3(B(x_0, R_0))}.
\end{align*}
Taking this inequality to the $ 3^{ {\rm rd}}$  power and integrate over $(-R_0^{5/2},0)$, we arrive at 
\[
\bigg(\intl_{ Q((x_0,0), R_0)} | \nabla p_h|^3 \phi ^3   dxd\tau\bigg)^{ 1/3} \le c R_2^{ -1} 
\bigg(\intl_{ Q((x_0,0), R_0)} | v|^3 \phi ^3   dxd\tau\bigg)^{ 1/3},
\] 
which shows that 
\[
II \le R^{ -1}_2  \intl_{ Q((x_0,0), R_0)} | v|^3   dxd\tau \le c(\delta \var )^{ 3/2}. 
\]
Inserting the estimates of $ I$ and $ II$ into \eqref{locen3}, we obtain 
 \begin{align}
&\intl_{ B(x_0, R_0)} | \widetilde{v} (t)|^2 \phi dx \le  | B(1)|^{ 1/3} \delta \var
+ c(\delta \var )^{ 3/2}. 
\label{locen33}
\end{align}   

Let $ -\frac12 < \tau <0$ be specified below. Once more using the fact that $ \nabla p_h$ is harmonic, we find 
for all $ \tau \in (-R_0^{5/2},0)$
\begin{align}
\| \nabla p_h(\tau )\|^2_{ L^2(B(x_0, \tau R_0))} &\le c (-\tau R_0)^{ 3} \| \nabla p_h(\tau )\|^2_{ L^\infty(B(x_0,  R_0/2))}
\cr
&\le c (-\tau )^{ 3} \| \nabla p_h(\tau )\|^2_{ L^2(B(x_0, R_0))} 
\cr
&\le c (-\tau )^{ 3} \| v(\tau )\|_{ L^2(B(1))}^2 = 
c (-\tau) ^{ 3}E,
\label{pressen}
\end{align}
where $ E= \| v\|^2_{ L^\infty(-1,0; L^2(B(1)))}$. Combining \eqref{locen33} and \eqref{pressen} with the choice $ \tau = -(\delta \var)^{ 1/3} $, we deduce that 
 \begin{align}
&\intl_{ B(x_0, -\tau R_0)} | v (t)|^2 \phi dx \le  c(| B(1)|^{ 1/3} + E+1) \delta \var. 
\label{locen4}
\end{align}   
In the above estimate we may choose $ \delta = \frac{1}{c(| B(1)|^{ 1/3} + E+1)}$ such that the desired estimate follows 
with $ R_1 = (\delta \var)^{ 1/3} R_0$.  \hfill \Beweisende

 \vspace{0.5cm}  
{\bf Proof of Theorem\,\ref{thm2.1}}:   Let $ \var > 0$ be arbitrarily chosen. 
By virtue of Lemma\,\ref{lem2.2} for every $ x_0\in B(1/2)$ 
there exists $ 0< R_1(x_0)= R_1(\var , x_0)< 1$ such that 
\[
\intl_{B(x_0, R_1(x_0))} | v (t)|^2  dx \le  \var  \quad  \forall\, -R_1(x_0)^{ 5/2} \le t < 0.  
\]    
Since $ \overline{B(1/2)}$ is compact,  we find a finite sequence of points $ \{ x_1, \ldots, x_m\}$ such that 
$ \{ B(x_i, R_1(x_i) /2)\}$ covers $\overline{B(1/2)} $.  Let $ y_0\in B(1/2)$ be arbitrarily chosen. There exists 
$ i\in \{ 1, \ldots, m\}$ with $ y_0 \in  B(x_i, R_1(x_i) /2)$. Obviously, $  B(y_0, R_1(x_i) /2) \subset  B(x_i, R_1(x_i) )$ and 
thus, 
\[
\intl_{B(y_0, R_1(x_i)/2)} | v (t)|^2  dx \le  \var  \quad  \forall\, -R_1(x_i)^{ 5/2} \le t < 0.  
\]    
Setting $ \widetilde{R} = \frac{1}{2} \min \{ R_1(x_1), \ldots, R_1(x_m)\}$,  we deduce that for all  $ y_0\in \overline{B(1/2)} $
it holds 
\begin{equation}
\intl_{B(y_0, \widetilde{R}) } | v (t)|^2  dx \le  \var\quad   \quad  \forall\, -\widetilde{R} ^{ 5/2} \le t < 0.  
\label{energy5}
\end{equation}
This completes the proof of assertion of the theorem.  \hfill \Beweisende 
 
 \vspace{0.5cm}  
 As an immediate consequence of Theorem\,\ref{thm2.1} we get the following smallness result for the $ L^\infty$ blow-up.
 \begin{cor}
 \label{cor2.5}
 Let the assumptions of Theorem\,\ref{thm2.1} be satisfied. Then for every $ \var >0$ there exists $ t_0=t_0(\var )\in(-1,0) $ such that 
 \begin{equation}
 \| v(t)\|_{ L^\infty(B(1/2))} \le \var (-t)^{ -3/5}\quad  \forall\,t_0 \le t < 0. 
 \label{epsLi}
 \end{equation}
  \end{cor}
 
 {\bf Proof}: Let $ \var >0 $ be arbitrarily chosen. Let $ 0<\delta<1 $ be fixed, which will be specified below.  We apply Theorem\,\ref{thm2.1} 
 with $\var _0=( \delta \var )^5$ in place of $ \var $. Let $ \widetilde{R} = \widetilde{R} (\var _0)$ such that \eqref{morrey} holds true 
 for $ \var _0$ in place of $ \var $. Applying the Gagliardo-Nirenberg inequality \eqref{A.1} with $ n=3, p=2$ and $ q=\infty$ together 
 with \eqref{morrey},   we obtain  
 \begin{align}
\| v(t)\|_{ L^\infty(B(\widetilde{R} ))} &\le c  \Big\{\widetilde{R} ^{-3/2 }\| v(t)\|_{ L^2(B(\widetilde{R} ))} + 
\| v(t)\|_{ L^2(B(\widetilde{R} ))}^{ 2/5}
 \| \nabla v(t)\|_{ L^{\infty}(B(\widetilde{R} ))}^{ 3/5}\Big\}. 
 \cr
&\le c \widetilde{R} ^{ -3/2} E^{ 1/2} + c C_0\delta \var (-t)^{ -3/5}   
\cr
&\le c \left\{\widetilde{R} ^{ -3/2} E^{ 1/2} (-t)^{ 3/5} +  C_0\delta \var \right\}(-t)^{ -3/5}.  
 \label{2.20}
 \end{align}
 We may choose $ \delta = \frac{1}{2 cC_0}$, and then $ t_0\in (-1,0)$ so that 
 $ \widetilde{R} ^{ -3/2} E^{ 1/2} (-t_0)^{ 3/5} \le \frac{\var }{2}$.  Then, \eqref{epsLi} follows from \eqref{2.20}.  \hfill \Beweisende  
  
  \vspace{0.3cm}
Using the  Gagliardo-Nirenberg inequality \eqref{GaNi6} instead of  \eqref{A.1} in the proof  Corolloary \ref{cor2.5}, 
we also get the uniform smallness  of the H\"older norm for any H\"older exponent $\gamma \in (0,1)$. 

\begin{cor}
\label{cor2.6}
Let the assumptions of Theorem\,\ref{thm2.1} be satisfied. Then for every $ 0< \gamma < 1$ and  for every $ \var >0$ there exists 
$ t_1=t_1(\gamma , \var )\in(-1,0) $ such that 
 \begin{equation}
 [v(t)]_{ C^{ 0, \gamma }} \le \var (-t)^{ -1 + \frac{2}{5}(1-\gamma )}\quad  \forall\,t_1 \le t < 0. 
 \label{epsHo}
 \end{equation}

\end{cor}

 \section{Local estimate for the second gradient}
 \label{sec:-4}
 \setcounter{secnum}{\value{section} \setcounter{equation}{0}
 \renewcommand{\theequation}{\mbox{\arabic{secnum}.\arabic{equation}}}}
 
 \begin{lem}
 \label{lem4.1}
 Let $ 2 \leq p< +\infty$ and $ v\in L^\infty(-1, 0; L^2(B(1)))\cap L^\infty(-1,0[); W^{1,\, \infty}(B(1)))$ be a solution to the Euler equations satisfying
 \begin{equation}
\sup_{ t\in (-1,0)} (-t) \| \nabla v(t)\|_{ L^\infty(B(1))} \le C_0<+\infty, 
 \label{4.1}
 \end{equation} 
together with $ v(-1)\in W^{2,\, p}(B(1))$, then  we have
 \begin{equation}
 \sup _{ t\in (-1,0)}(-t)^{ cC_0 p^2} \| \nabla ^2 v(t)\|_{ L^p(B(1/8))} < +\infty,
 \label{4.1a}
 \end{equation}
 where $c>0$ is an absolute constant.
 \end{lem}

 {\bf Proof}:  According to  Corollary \ref{cor2.5} we can assume without the loss of generality that 
 \begin{equation}
 \sup_{ t\in (-1,0)} (-t)^{ 3/5} \| v(t)\|_{ L^\infty(B(1))} \le \frac{1}{32}. 
 \label{4.4b}
 \end{equation} 
We set  $ \theta = \frac{1}{4}$, and  define 
 \[
w(y, t) := v(y+(-t)^\theta y, t ),\quad  (y,t)\in \R^{3}\times (-1,0).  
\]
 Clearly, $ w$ solves the  following modified Euler equations in $ \R^{3}\times (-1,0)$.
  \begin{align}
  &\partial _t w + \frac{\theta  (-t)^{ \theta -1}y}{1+ (-t)^\theta } \cdot \nabla w + \frac{1}{1+(-t)^\theta } (w\cdot \nabla) w
 = -\nabla \pi ,
 \label{4.3}\\
 &\qquad \nabla \cdot w =0.
 \label{4.2}
 \end{align}
Let us set  $ \Omega =\nabla \times w$. Then, we find from \eqref{4.3} that  $ \Omega $ solves the equation. 
 \begin{align}
 &\partial _t \Omega  + \frac{\theta (-t)^{ \theta -1}}{1+ (-t)^\theta } \Omega  + 
 \frac{\theta (-t)^{ \theta -1}y}{1+ (-t)^\theta } \cdot \nabla \Omega  + \frac{1}{1+(-t)^\theta } (w\cdot \nabla) \Omega 
 \cr
 &\qquad \qquad =  \frac{1}{1+(-t)^\theta }  \Omega \cdot \nabla w\quad  \text{ in}\quad\R^{3}\times (-1,0). 
 \label{4.4}
 \end{align}
Applying  the derivative $ \partial _i$, $ i=1,2,3$,  to the both sides of \eqref{4.4}, we see that 
$$ U_i=\partial _i \Omega = \nabla\times \partial _i w$$  solves the equations 
 \begin{align}
 &\partial _t U_i  + 2\frac{\theta (-t)^{ \theta -1}}{1+ (-t)^\theta } U_i  + 
 \frac{\theta (-t)^{ \theta -1}y}{1+ (-t)^\theta } \cdot \nabla U_i  + \frac{1}{1+(-t)^\theta } (w\cdot \nabla) U_i
 \cr
 &\qquad \qquad = \frac{1}{1+(-t)^\theta } \Big\{-(\partial _i w\cdot \nabla) \Omega  +   U_i \cdot \nabla w
 + \Omega \cdot \nabla \partial _i w\Big\}
 \cr
 &\quad  \text{ in}\quad\R^{3}\times (-1,0). 
 \label{4.5}
 \end{align}
Note that $ x:= (1+ (-t)^\theta ) y\in B(1)$ for $ y \in B(1/2)$.
 Let $  \eta \in C^{\infty}([0, +\infty))$ be non increasing such that 
 $ \eta \equiv 1$ on $ \Big[0, \frac{1}{4})$, $\eta \equiv 0 $ on  $[\frac12, +\infty)$. We set $ \phi (y)= \eta (| y|)$.  Let 
 $ 2 \le p <+\infty$. 
 We multiply \eqref{4.5} by $ U_i | U|^{ p-2} \phi $, taking the sum from $ i=1$ to $ 3$,  integrating it  over $ B(1/2)\times (-1, t)$ with $ t\in (-1, 0)$,  and applying the integration by parts, we  obtain
 \begin{align}
 &\intl_{B(1/2)} | U(t)|^p \phi dy+ (2p- 3)  \intl_{-1}^{t}  \intl_{B(1/2)} \frac{\theta (-s)^{ \theta -1}}{1+ (-s)^\theta } | U|^p \phi    dy   dt 
   \cr
 &  -    \intl_{-1}^{t}  \intl_{B(1/2)  \setminus B(1/4)} \left\{ \frac{\theta (-s)^{ \theta -1}}{1+ (-s)^\theta } | U|^p   \eta ' (| y|) | y|    + 
 \frac{1}{1+ (-s)^\theta }
 | U|^p  \frac{w\cdot y}{| y|} \eta'(| y|)\right\} dy   d s
  \cr
 &  =  \intl_{B(1/2)} | U(-1)|^p \phi dy 
 \cr
  &\qquad +   p \intl_{-1}^{t}  \intl_{B(1/2)}
 \frac{1}{1+(-s)^\theta }\Big\{- (\partial _i w\cdot \nabla) \Omega  +   U_i \cdot \nabla w
 + \Omega \cdot \nabla \partial _i w\Big\}\cdot U_i | U|^{ p-2}\phi dy   ds. 
 \label{4.6}
 \end{align}   
 
In order to get the positive sign for the third  term on the left-hand side of \eqref{4.6} we use  \eqref{4.4b},  which implies that 
for all $ y\in B(1/2)  \setminus B(1/4)$ it holds 
\begin{align*}
\theta (-s)^{ \theta -1} | y|    +    \frac{w\cdot y}{| y|} 
\ge \frac{1}{16} (-s)^{ - 3/4} -  \frac{1}{32}(-s)^{ -3/5} \ge    \frac{1}{32} (-s)^{ - 3/4}. 
\end{align*}
Since $\eta' (|y|)\leq 0$ for all $y\in \Bbb R^3$, from \eqref{4.6}  we deduce the estimate 
 \begin{align*}
 &\intl_{B(1/2)} | U(t)|^p \phi dy
 \\
 & \le   \intl_{B(1/2)} | U(-1)|^p \phi dy  +  p  \intl_{-1}^{t}  \intl_{B(1/2)}
\frac{| \nabla w|}{1+ (-s)^\theta } (2|U|  + | \Omega| | \nabla ^2w| ) | U|^{ p-1} \phi  dy   ds. 
 \end{align*}   
 Observing the Type I condition for $ v$, applying H\"older's inequality and Young's inequality, and replacing $ \phi $ by $ \phi  ^p$, we get from the above inequality 
 \begin{align}
 &\intl_{B(1/2)} | U (t)|^p \phi^p dy
 \cr
 & \le   \intl_{B(1/2)} | U (-1)|^p \phi^p dy  +    pC_0  
 \intl_{-1}^{t}  \intl_{B(1/2)} (-s)^{ -1} (2| U(s) |^p + | \nabla ^2 w|| U(s) |^{ p-1} )\phi ^p dy   ds
 \cr
 & \le   \intl_{B(1/2)} |U (-1)|^p \phi^p dy  +  2 p C_0 
 \intl_{-1}^{t}  \intl_{B(1/2)} (-s)^{ -1} | U |^p  \phi ^p dy   ds
 \cr
 &\qquad \qquad +  pC_0\intl_{-1}^{t} (-s)^{ -1} \| \nabla^ 2 w(s) \phi \|_{ L^p(B(1/2)) } \| U(s) \phi \| ^{ p-1}_{ L^p(B(1/2)) }ds. 
 \label{4.7}
 \end{align}   
  Furthermore, using the Biot-Savart law and Calder\'on-Zygmund inequality  together with Lemma\,\ref{lemvort}, we get the   estimate    \begin{align}
&\| \nabla^ 2 w(s) \phi \|_{ L^p(B(1/2)) }
\cr
 &\quad \leq  c\sum_{i=1}^3  \| ( \nabla \partial_i  w (s))\phi\|_{ L^p(B(1/2)) } \cr 
&\quad \le  cC_p\sum_{i=1}^3 \| \partial_i  (\nabla \times w)(s) \phi \|_{ L^p(B(1/2)) }+ cC_p\|  \nabla w (s)\nabla \phi \|_{ L^p(B(1/2)) }
\cr 
& \quad \le c p \bigg(\intl_{B(1/2)} | U(s) |^p \phi^p dy\bigg)^{ 1/p} +c p  C_0 (-s)^{ -1}.  
\label{d2w}\end{align}
Hence, \eqref{4.7} together with Young's inequality yields  
\begin{align}
 &\intl_{B(1/2)} |U(t)|^p \phi^p dy
 \cr
 & \le  c C_0  p^2 \intl_{-1}^{t}  \intl_{B(1/2)} (-s)^{ -1}| U |^p \phi^p  dyds+ c C_0 p^2 (-t)^{-p} +   \intl_{B(1/2)} | U (-1)|^p \phi^p dy .  
 \label{4.8}
 \end{align}    
 We  define 
 \[
X(t) = \intl_{-1}^{t}  \intl_{B(1/2)} (-s)^{ -1}|U |^p \phi^p  dyds+\frac{ cC_0 p(-t)^{-p}}{ cC_0 -1} + \frac{1}{c C_0 p^2}  \intl_{B(1/2)} |  U (-1)|^p \phi^p dy.
\]
Thanks to \eqref{4.8} we find   
  \begin{align*}
X'(t) &= (-t)^{ -1}\intl_{B(1/2)} | U(t)|^p \phi^p dy + \frac{cC_0 p^2 }{cC_0-1} (-t)^{ -p-1}
  \\
&\le    c C_0 p^2(-t)^{ -1} \intl_{-1}^{t}  \intl_{B(1/2)} (-s)^{ -1}| U |^p \phi^p  dyds
\\
&\qquad  +  
cC_0 p^2(-t)^{ -p-1} +  (-t)^{ -1} \intl_{B(1/2)} | U (-1)|^p \phi^p dy  +  \frac{cC_0 p^2}{cC_0-1} (-t)^{ -p-1}
\\
&\leq cC_0 p^2  (-t)^{ -1}  X(t). 
\end{align*}
This shows that $t\mapsto  X(t)(-t)^{ c C_0 p^2}$ is non increasing on $(-1,0)$. Consequently,
\begin{align}
X(t) &\le (-t)^{ -cC_0 p^2} X(-1) 
\cr
&=  (-t)^{ -cC_0p^2}\bigg\{ \frac{cC_0p}{cC_0-1}  + \frac{1}{cC_0p^2}
\intl_{B(1/2)} | U(-1)|^p \phi^p dy\bigg\}. 
\label{4.9}
\end{align}
Combining \eqref{4.8} with  \eqref{4.9}, we arrive at 
\begin{align}
\| \nabla \Omega(t) \|_{ L^p(B(1/4))}^p  &\le cC_0 p^2  X(t) 
\cr
&\le  (-t)^{ -cC_0 p^2 } \bigg\{ \frac{c C_0^2 p^3}{ cC_0-1}  +
\| \nabla \Omega(-1) \|_{ L^p(B(1/2))}^p\bigg\}. 
\label{4.10}
\end{align}
 Once more, applying Lemma\,\ref{lemvort}, we get the assertion of the lemma.  \hfill \Beweisende

\vspace{0.3cm}
\hspace{0.5cm}
If we replace  the condition \eqref{4.1} by   
\begin{equation}
(-t)^\beta \| \nabla v(t)\|_{ L^\infty(B(1))} < +\infty, 
\label{4.11}
\end{equation}
for some $ 1< \beta < 1$, then  we  get the following bound for the  second gradient. 

 \begin{lem}
 \label{lem4.2}
 Let $ v\in L^\infty(-1, 0; L^2(B(1)))\cap L^\infty(-1,0[); W^{1,\, \infty}(B(1)))$ be a solution to the Euler equations with 
 $ v(-1)\in W^{2,\, p}(B(1))$, $ 3< p< +\infty$. We assume \eqref{4.11} holds for some $ 0< \beta <1$.   Then, 
 \begin{equation}
\| \nabla ^2 v(t)\|_{ L^\infty(-1,0; L^p(B(1/8))} < +\infty.   
 \label{4.12}
 \end{equation}
 \end{lem}

{\bf   Proof}: Repeating the proof of  Lemma\,\ref{lem4.1} up to \eqref{4.7}, and using \eqref{4.11} instead of  \eqref{4.1}, we 
obtain 
\begin{align}
 &\intl_{B(1/2)} | U (t)|^p \phi^p dy
 \cr
  & \le   \intl_{B(1/2)} |U (-1)|^p \phi^p dy  +  c
 \intl_{-1}^{t}  \intl_{B(1/2)} (-s)^{ -\beta } | U |^p  \phi ^p dy   ds
 \cr
 &\qquad \qquad +  c\intl_{-1}^{t}  (-s)^{ -\beta } \| \nabla^ 2 w(s) \phi \|_{ L^p(B(1/2))} \| U(s) \phi \| ^{ p-1}_{ L^p(B(1/2))}ds. 
 \label{4.13}
 \end{align}   
 As in the proof of Lemma\,\ref{lem4.1}, using  \eqref{d2w} Lemma\,\ref{lemvort} and Young's inequality, we are led to  
   \begin{align}
 &\intl_{B(1/2)} |U(t)|^p \phi^p dy
 \cr
 & \le  c\intl_{-1}^{t}  \intl_{B(1/2)} (-s)^{ -\beta }| U |^p \phi^p dyds +  c (-t)^{-p -\beta+1}+    \intl_{B(1/2)} | U (-1)|^p \phi^p dy.  
 \label{4.14}
 \end{align}    
 We  define 
 \[
\hspace*{-1.5cm}X(t) = \intl_{-1}^{t}  \intl_{B(1/2)} (-s)^{ -\beta }|U |^p \phi^p  dyds+  \frac{c}{(p+\beta-1)(c-1)} (-t)^{ -p-\beta+1 } +  \intl_{B(1/2)} |  U (-1)|^p \phi^p dy.
\]
In view of  \eqref{4.14} we obtain   
  \begin{align*}
X'(t) &= c (p+\beta-1)(-t)^{ -\beta }X(t)\leq c (-t)^{ -\beta }X(t). 
\end{align*}
By means of Gronwall's lemma we find 
\begin{align}
X(t) &\le X(-1)e ^{ c\int_{-1} ^t (-s)^{-\beta} ds }.  
\label{4.15}
\end{align}
Combining \eqref{4.15} with  \eqref{4.14}, we arrive at 
\begin{align}
\label{domega}
\| \nabla \Omega(t) \|_{ L^p(B(1/4))}^p  &\le c X(t) \le X(-1)e ^{ c \intl_{-1}^{t} (-s)^{ -\beta }ds}  
\\
&\le \Big(c + \| \nabla \Omega(-1) \|_{ L^p(B(1/4))}^p\Big)e ^{ \frac{c}{1-\beta }}. 
\end{align}
Applying Lemma\,\ref{lemvort}, we get the assertion of the lemma.  \hfill \Beweisende 

  \section{Proof of Theorem\,\ref{thm1.1}}
  \label{sec:-6}
  \setcounter{secnum}{\value{section} \setcounter{equation}{0}
  \renewcommand{\theequation}{\mbox{\arabic{secnum}.\arabic{equation}}}}
 
 The hypothesis \eqref{decay} implies that there exists $\eta \in (0, r_0)$ such that 
 \begin{equation} 
 \label{hypo}\sup_{-\eta<t< 0} (-t)\|\nabla v(t)\|_{L^\infty ( B(\eta ))} < 1. 
 \end{equation}  
 Then, by rescaling, one may assume without the loss of generality that $\eta =1$ in \eqref{hypo}. 
 Let $ 0< \rho < r_0$ be fixed.  Hence, it will sufficient  to show that 
for every $ x_0\in \overline{B(\rho )} $  it holds   $\nabla v \in L^\infty (-1,0; W^{2,\, p_0}(B(r/32)))$, where $ r= r_0-\rho $.  
We define the rescaled velocity by means of 
\[
\widetilde{v} (y,t) = r^{ 3/2} v(x_0+ ry,  r^{ 5/2} t),\quad  (y,t) \in B(1)\times (-1,0).  
\]
For notational simplicity we write again $ v$ in place of $ \widetilde{v} $ and prove  that 
$$ \nabla v \in 
L^\infty (-1,0; W^{2,\, p_0}(B(1/32))).$$

\hspace{0.5cm}
Thanks to Corollary\,\ref{cor2.5} we may assume that 
\begin{equation}
\sup_{ t\in (-1,0)} (-t)^{ -3/5} \| v(t)\|_{ L^\infty(B(1/2))} <\frac{1}{10}.  
\label{blowupLi}
\end{equation}
Let $ t_0\in (-1,0)$ be arbitrarily chosen but fixed. Let $ x_0\in B(1/4)$. By $ X(x_0 , t_0;  s)$ we denote the trajectory of the particle 
which is located at $ x_0$ at time $ s=t_0 $. More precisely, $ s \mapsto X(x_0 , t_0;  s)$ solves the following ODE
\begin{equation}
\frac{d X}{d s}(x_0, t_0;  s)  = v(X(x_0 , t_0;  s), s)\quad  \text{ in }\quad  [-1,0),\quad  X(x_0 , t_0;  t_0) = x_0.  
\label{6.1}
\end{equation}
Since $ v(t)$ is Lipschitz in $ B(1)$ for all $ t\in (-1,0)$ we first get a local solution of \eqref{6.1} in some maximal  interval 
$I=(t_1, t_2)$, such that $ X(x_0 , t_0;  s) \in B(1/2)$ for all $ s\in I$.  We claim that $ I = (-1,0)$.  In fact integration 
over $ (t, t_0)$  of \eqref{6.1} for some $ t\in I$ yields 
\[
X(x_0 , t_0;  t) - x_0 =  \intl_{t_0}^{t} v(X(x_0 , t_0;  s), s)  ds. 
\]
Using  the triangle inequality along with \eqref{blowupLi}, we obtain 
\begin{equation}
| X(x_0 , t_0;  t) - x_0| <\frac{1}{10} \intl_{t}^{t_0} (-s)^{-3/5}  ds \le \frac{1}{4}. 
\label{control}
\end{equation} 
Thus  $ I \neq (-1,0)$ would lead to a contradiction,  since by  \eqref{control} we may extend the solution to a larger interval,  which 
violates the maximal property of $ I$. This shows that the whole trajectory $ X(x_0 , t_0;  t) - x_0$ remains in $ B(1/2)$ for all $ t\in (-1,0)$. 

Let $ \omega = \nabla \times v$, which  solves the vorticity equations
\begin{equation}
\partial _t \omega + v\cdot \nabla \omega = \omega \cdot \nabla v\quad  \text{ in}\quad  \R^{3}\times (-1,0). 
\label{6.2}
\end{equation}
Observing \eqref{6.1}, by means of the chain rule we infer from \eqref{6.2}  that $s \mapsto \omega (X(x_0, t_0;  s), s)$ solves the following ordinary differential equation
\[
\frac{d }{d s} \omega (X(x_0, t_0;  s), s) = \omega(X(x_0, t_0;  s), s) \cdot \nabla v(X(x_0, t_0;  s), s) \quad  \text{ in}\quad  (-1,0).
\]
Multiplying the above equation by $ \omega (X(x_0, t_0;  s), s)$,  we see that the function $\psi (s):= | \omega (X(x_0, t_0;  s), s)|^2 $ 
satisfies the inequality 
\[
\frac{1}{2}\psi ' \le  | \nabla v(X(x_0, t_0;  s), s)|\psi\quad  \text{ in}\quad  (-1,0), 
\]
 Observing the assumption \eqref{decay}, the above inequality implies for some $ 0<\delta< \frac{2}{5} $
 \begin{equation}
 \psi ' \le 2 (-s)^{ -1} (1-\delta ) \psi\quad   \text{ in}\quad  (-1,0).  
 \label{6.5}
 \end{equation} 
 This, immediately shows that $ (-s)^{ -2(1-\delta )}\psi (s)$ is non increasing in $ (-1,0)$. Accordingly, 
 \begin{equation}
\psi (s) \le    (-s)^{ -2(1-\delta )} \psi (-1)\quad  \forall\, s \in (-1,0). 
 \label{6.6}
 \end{equation}
 In particular, inequality \eqref{6.6} with $ s=t_0$ yields  \[
 | \omega (x_0, t_0)| \le (-t_0)^{ -1+\delta } | \omega (X(x_0, t_0;  -1), -1)| \le  (-t_0)^{ -1+\delta } \| \omega (-1)\|_{ L^\infty(B(1))}. 
\]
 This estimate gives 
 \begin{equation}
\| \omega (t)\|_{ L^\infty(B(1/4))} \le  (-t)^{ -1+\delta } \| \omega (-1)\|_{ L^\infty(B(1))}\quad  \forall\,t\in (-1,0). 
 \label{6.7}
 \end{equation} 
 Applying \eqref{vort},  we infer from \eqref{6.7} together with \eqref{blowupLi}  that the following  estimate holds  
 for all $ 3< q< +\infty$, 
 \begin{equation}
\| \nabla v (t)\|_{ L^q(B(1/8))} \le  cq \Big\{(-t)^{ -1+\delta } \| \omega (-1)\|_{ L^\infty(B(1))} 
+  (-t)^{ -3/5}\Big\}\quad \forall\,t\in (-1,0). 
 \label{6.8}
 \end{equation}    
 Noting  that $ \frac{3}{5} < 1-\delta $, by means of Sobolev's embedding theorem we deduce form \eqref{6.8} for 
 $ \gamma = 1- \frac{3}{q}$
   \begin{equation}
[v(t)]_{ C^{ 0, \gamma }(B(1/8))} \le  c (1-\gamma )^{ -1} (-t)^{ -1+\delta } \Big(1+ \| \omega (-1)\|_{ L^\infty(B(1))} \Big)
\quad \forall\,t\in (-1,0)
 \label{6.9}
 \end{equation}   
 with a constant $ c>0$, which  remains bounded as $ q \rightarrow +\infty$.  
 
 \hspace{0.5cm}
 Appealing to \eqref{GaNihol} (cf. Lemma\,A.2) with $ n=3$, we see that for all  $ t\in (-1,0)$ the following inequality holds true
 \begin{align*}
\| \nabla v(t)\|_{ L^\infty(B(1/8))} &\le  c  [v(t)]_{ C^{ 0, \gamma }(B(1/8))} 
  \\
 &\qquad + c 
\Big([v(t)]_{ C^{ 0, \gamma }(B(1/8))}\Big) ^{ 1- \frac{1-\gamma }{2-\gamma - 3/p_0}}
\|\nabla^2 v(t) \|_{ L^q(B(1/8))}^{ \frac{1-\gamma }{2-\gamma - 3/p_0}}. 
 \end{align*}
We estimate the right-hand side of the above inequality by the aid of   \eqref{6.9} and \eqref{4.1a}. This gives 
\begin{align*}
&\| \nabla v(t)\|_{ L^\infty(B(1/8))} 
\\
&\quad \le c (1-\gamma )^{ -1}(-t)^{ -1+\delta } \Big(1+ \| \omega (-1)\|_{ L^\infty(B(1))} \Big)
\\
&\qquad + c (1- \gamma )^{  1- \frac{1-\gamma }{2-\gamma - 3/p_0}} (-t)^{ -1+\delta + (-C_2+1-\delta )\frac{1-\gamma }{2-\gamma - 3/p_0} } 
\Big(1+ \| \omega (-1)\|_{ L^\infty(B(1))} \Big) ^{ 1- \frac{1-\gamma }{2-\gamma - 3/p_0}},
\end{align*}    
where $ C_2= c C_0p_0^2$ stands for the constant in Lemma\,\ref{lem4.1}.  We can choose  $ \gamma \in (0,1)$ such that 
 \[
 (-C_2+1-\delta )\frac{1-\gamma }{2-\gamma - 3/p_0} \geq - \frac{\delta }{2}.
\]
  With this choice of $ \gamma $ we get 
  \begin{equation}
  (-t)^{ 1- \frac{\delta}{2}} \| \nabla v(t)\|_{ L^\infty(B(1/8))} < +\infty. 
  \label{6.10}
    \end{equation}
  Thus, we are in a position to apply Lemma\,\ref{lem4.2} for $ \beta = 1- \frac{\delta }{2}$,  which yields 
  \[
v \in L^\infty(-1,0; W^{2,\, p_0}(B(1/32))). 
\]
This completes the proof of the theorem. 
 \hfill \Beweisende 
  
 \section{Proof of Theorem\,\ref{thm1.2}}
 \label{sec:-?}
 \setcounter{secnum}{\value{section} \setcounter{equation}{0}
 \renewcommand{\theequation}{\mbox{\arabic{secnum}.\arabic{equation}}}}

Let  $ \zeta \in C^{\infty}_c(B(1))$ denote a cut off function such that $ 0 \le \zeta  \le 1$ in $ B(1)$, 
$ \zeta \equiv 1 $ on $ B(1/2)$, and $ | \nabla \zeta | \le c $.  For $ 0< r< +\infty$, and $ x_0\in \R^{3}$ we define  
\[
\zeta _r= \zeta _{r} (x) = \zeta \Big(\frac{x-x_0}{r}\Big),\quad  x\in \R^{3}.
\]  
 We set $ R:= \frac{r_0}{2}$. Clearly, $ \zeta _{ R} \in C^{\infty}_c(B(x_0, R))$ is a cut off function  on $ B(R)= B(x_0, R)$ with $ | \nabla \zeta _{ B(R)}| \le c R^{ -1}$. 
 We define the modified mean value 
 \[
\widetilde{f}_{ B(R)} = \frac{1}{ \intl_{B(R)} \zeta_R dx } \intl_{B(R)} f \zeta_R dx.   
\]
Let $ t\in (-1, 0)$, and $ 0< R \le  r_0 (-t)^{ 2/5}$ be fixed. For $ x_0\in B(R)$ we get 
\begin{align*}
| \nabla v(x_0, t)| &\le  | \nabla v(x_0, t) -  \widetilde{\nabla v( t)}_{ B(R)} | + | \widetilde{\nabla v( t)}_{ B(R)}| 
\\
&\le \frac{1}{ \intl_{B(R)} \zeta_R dx } \intl_{B(R)} | \nabla v(x_0, t) -\nabla v(x,t)| \zeta_R  dx 
\\
&\qquad \qquad \qquad \qquad \qquad + \frac{1}{ \intl_{B(R)} \zeta_R dx}
\bigg| \intl_{B(R)} \nabla v(x,t)\zeta_R dx\bigg|
\\
&\le \osc_{ B(x_0, R)} (\nabla v(t)) + \frac{1}{ \intl_{B(R)} \zeta_R dx}
\bigg| \intl_{B(R)} \nabla v(x,t)\zeta_R dx\bigg|. 
\end{align*}
Applying the integration by parts and H\"{o}lder's  inequality,   we find 
\begin{align*}
\bigg| \intl_{B(R)} \partial _i v(x,t)\zeta_R dx\bigg| &= \bigg| \intl_{B(R)}  v(x,t) \partial _i\zeta_R dx\bigg|
\le c R^{ 1/2} \| v(t)\|_{ L^2(B(R))}.  
\end{align*}
This leads to the inequality 
\begin{equation}
| \nabla v(x_0, t)| \le  \osc_{ B(x_0, R)} (\nabla v(t)) + c R^{ -5/2} \| v(t)\|_{ L^2(B(R))}. 
\label{3.1}
\end{equation}
In particular, observing \eqref{decay} from \eqref{3.1} with  $ R= r_0 (-t)^{ 2/5}$,  we deduce 
\begin{equation}
(-t)| \nabla v(x_0, t)| \le 1-\delta  + c r_0^{ -5/2}\| v(t)\|_{ L^2(B(x_0, r_0 (-t)^{ 2/5}))}. 
\label{3.2}
\end{equation}
Since $\| v(t)\|_{ L^2(B(r_0 (-t)^{ 2/5}))} $ is bounded by $ E^{ 1/2} = \| v\|_{ L^\infty(-1,0;L^2(B(1))}$,  we immediately 
get from \eqref{3.2}   that $ v$ has Type I blow up at $ t=0$ with respect to the velocity gradient. This allows us to apply Theorem\,\ref{thm2.1} which  yields the existence  of $ \widetilde{R} $ such that for all $ t\in (-\widetilde{R}^{ 5/2},0 )$ and for all $ x_0\in B(1/2)$ it holds 
\[
c r_0^{ -5/2}\| v(t)\|_{ L^2(B(x_0, r_0 (-t)^{ 2/5}))} \le  c r_0^{ -5/2}\| v(t)\|_{ L^2(B(\widetilde{R}) )} \le \frac{\delta }{2},
\] 
which gives 
\begin{equation}
\sup_{ t\in (-\widetilde{R}^{ 5/2}, 0 )}  (-t)\| \nabla v(t)\|_{ L^\infty(B(r_0/2))} \le 1- \frac{\delta }{2}. 
\label{3.3}
\end{equation}
This shows that the condition \eqref{decay}  in Theorem\,\ref{thm1.1} is satisfied,  which  yields that 
$  v\in L^\infty(-1,0; W^{2,\, p_0}(B( r_0/4)))$. This completes the proof of the theorem.  \hfill \Beweisende

$$\mbox{\bf Acknowledgements}$$
Chae was partially supported by NRF grants 2016R1A2B3011647, while Wolf has been supported 
supported by the German Research Foundation (DFG) through the project WO1988/1-1; 612414.  
  
\appendix
\section{Gagliardo-Nirenberg's inequality on a ball}
\label{sec:-A}
\setcounter{secnum}{\value{section} \setcounter{equation}{0}
\renewcommand{\theequation}{\mbox{A.\arabic{equation}}}}

\begin{lem}
\label{lemGaNi1}
Let $ 1 \le p, q < +\infty$. We assume $ q >n$. Let $ B(R)= B(x_0, R)$ be any ball.  
Then for all $ f\in L^p (B(R))\cap W^{1,\, q}(B(R))$ 
it holds
\begin{equation}
\| f\|_{ L^\infty(B(R))} \le c R^{ -n/p}\| f\|_{ L^p(B(R))} + c \| f\|_{ L^p(B(R))}^{ 1- { \frac{nq}{pq- pn+ qn}}} \| \nabla f\|_{ L^q(B(R))}^{ { \frac{nq}{pq- pn+ qn}}},
\label{A.1}
\end{equation}
where the constant in \eqref{A.1} depends only on $ p,q$ and $ n$ but not on $ R>0$ and $ x_0$.
\end{lem}

{\bf Proof}: Let $ 0 < \lambda  \le 1$.  By means of Sobolev's embedding theorem we find for all $ 0< \lambda  \le 1$ 
\begin{align}
\label{sobolev}
\| f\|_{ L^\infty(B(1)\cap B(x_0, \lambda ))} \le c\lambda ^{ -n/p} \| f\|_{ L^p(B(1))} + 
c\lambda ^{ 1- n/q}\|\nabla f \|_{ L^q(B(1))}.
\end{align}
In case $ \|\nabla f \|_{ L^q(B(1 ))}  \le \| f\|_{ L^p(B(1))} $ the assertion is trivially fulfilled by setting $\lambda =1$ in \eqref{sobolev}.  In the opposite case we 
choose $ \lambda $ such that $ \lambda ^{ -n/p} \| f\|_{ L^p(B(\lambda ))} = \lambda ^{ 1- n/q}
\|\nabla f \|_{ L^q(B(\lambda ))}$, i.e. 
\[
\lambda  =\bigg( \frac{  \| f\|_{ L^p(B(1))}  }{  \|\nabla f \|_{ L^q(B(1 ))  }}\bigg)^{ \frac{1}{1-n/q+ n/p}} \le 1. 
\]
This implies that 
\[
\| f\|_{ L^\infty(B(1))} \le c \| f\|_{ L^p(B(1))} + c \| f\|_{ L^p(B(1))}^{ 1- { \frac{nq}{pq- pn+ qn}}} \| \nabla f\|_{ L^q(B(1))}^{ { \frac{nq}{pq- pn+ qn}}}.
\]

The assertion now follows easily by means of a standard scaling argument.  \hfill \Beweisende 

\vspace{0.3cm}
The following lemma provides an estimate between H\"older space and Sobolev space.

\begin{lem}
\label{lemGaNi2}  
Let $ 1 < p < +\infty$ with $ p > n$, and let $ 0< \gamma   < 1$. Then for every ball and 
$ f\in  C^{ 0,\gamma }(\overline{B(R)} )\cap W^{2,\, p}(B(R))$ it holds 
\begin{equation}
\| \nabla f\|_{ L^\infty(B(R))} \le  c R^{ \gamma -1} [f]_{ C^{ 0, \gamma }(B(R))} + c [f]_{ C^{ 0, \gamma }(B(R))} ^{ 1- \frac{1-\gamma }{2-\gamma - n/p}}
\|\nabla^2 f \|_{ L^q(B(R)}^{ \frac{1-\gamma }{2-\gamma - n/p}}. 
\label{GaNihol}
\end{equation}
\end{lem} 
 
{\bf Proof}: As in the proof of Lemma\,\ref{lemGaNi1} it will be sufficient to prove the assertion for the unite ball $ B(1)$. 
The general case easily follows by a standard  scaling argument.  By means of Sobolev's embedding theorem we find for all $ 0< \lambda  \le 1$ 
\begin{align*}
\lambda\| \nabla  f\|_{L^\infty(B(1)\cap B(x_0, \lambda ))} \le c\lambda ^{ \gamma }  [f]_{ C^{ 0, \gamma }(\overline{B(1)} } + 
c\lambda ^{ 2- n/q}\|\nabla^2 f \|_{ L^q(B(1))}. 
\end{align*}
In case $ \|\nabla^2 f \|_{ L^q(B(1))}\le  [f]_{ C^{ 0, \gamma }(\overline{B(1)} )}$ the assertion is trivially fulfilled. Otherwise, we set 
\[
\lambda  =\bigg( \frac{ [f]_{ C^{ 0, \gamma }(\overline{B(1)} )}}{\|\nabla^2 f \|_{ L^q(B(1))}}\bigg)^{ \frac{1}{2-\gamma - n/q}}. 
\]
Then the above inequality implies 
\[
\| \nabla f\|_{ L^\infty(B(1))} \le  c [f]_{ C^{ 0, \gamma }(B(1))} + c [f]_{ C^{ 0, \gamma }(B(1))} ^{ 1- \frac{1-\gamma }{2-\gamma - n/q}}
\|\nabla^2 f \|_{ L^q(B(1))}^{ \frac{1-\gamma }{2-\gamma - n/q}}. 
\]
This completes the proof of the lemma.  \hfill \Beweisende 

\vspace{0.5cm}  
Next, we provide the following elementary estimate of the H\"older semi norm. 

\begin{lem}
For all $ f\in L^\infty(\R^{n})\cap  W^{1,\, \infty}(\R^{n})$ it holds 
\begin{equation}
[f]_{ C^{ 0,\alpha }} \le 2^{ 1-\alpha }\| f\|_{ L^\infty} ^{ 1-\alpha } \| \nabla f\|^{ \alpha }_{ L^\infty}. 
\label{GaNi6}
\end{equation}
\end{lem}

{\bf Proof}: Elementary, for $ x,y\in \R^{n}$ with $ x \neq y$ we estimate 
\begin{align*}
\frac{| f(x)-f(y)|}{| x-y|^\alpha } &= | f(x)- f(y)|^{ 1-\alpha } \Big( \frac{| f(x)- f(y)|}{| x-y|}\Big)^\alpha 
\le 2^{ 1-\alpha } \| f\|_{ L^\infty} \| \nabla f\|_{ L^\infty}^\alpha. 
\end{align*}
After taking the supremum over all $ x, y\in \R^{n}$ with $ x \neq y$ on both sides of the above inequality, we get the assertion of the lemma.  \hfill \Beweisende

\vspace{0.3cm}
Using the well known Biot-Savart law together with Calder\'on-Zygmund's inequality\cite{ste},  we get the following localized inequality.

\begin{lem}
\label{lemvort}
Let $ 1< p< +\infty$. Then for every ball  $ B(R) \subset  \R^{3}$ and $ u\in W^{1,\, p}(B(R))$ it holds 
\begin{equation}
\| \nabla u \phi \|_{ L^p(B(R))} \le C_p \Big(\| \nabla \times u \phi \|_{ L^p(B(R))} + \| \nabla \phi \|_{ \infty} \| u\|_{ L^p(B(R))}\Big),
\label{vort}
\end{equation}
for all non negative $\phi \in C^{\infty}_{\rm c}(B(R)) $,  where $ C_p\leq c p$  for $p\in [2, \infty)$, while  $C_p \leq \frac{c}{p-1}$  for $p\in (1, 2]$ with  $c$  independent  of $p$.
\end{lem}

{\bf Proof}:  Let $ \phi \in C^{\infty}_{\rm c}(B(R))$ be a non negative function.  Since 
$
(\nabla   u) \phi =  \nabla  (u \phi) -u \otimes  \nabla  \phi, 
$
by the Biot-Savart law and Calder\'on-Zygmund inequality we get 
\begin{align*}
\|  \nabla   (u\phi) \|_{ L^p } &\le C_p  \|   \nabla \times  (u\phi) \|_{ p } 
\\
&=C_p  \|   \nabla \times  u \phi - u \times \nabla \phi  \|_{ L^p } 
\\
&\le 
C_p \|   \nabla \times  u \phi   \|_{ L^p } +  C_p  \| \nabla \phi \|_{ \infty} \| u  \|_{ L^p(B(R)) }.
\end{align*}
This immediately leads to \eqref{vort}.  \hfill \Beweisende  

\bibliographystyle{siam}

\end{document}